\ifx\shlhetal\undefinedcontrolsequence\let\shlhetal\relax\fi

\documentstyle[12pt,psfig]{article}

\input amssym.def
\input amssym.tex
\define\peg{{}^\frown\!}

\define\cf{\; {\rm cf} \; }
\define\ife{\underbar{if} }
\define\then{\underbar{then} }
\define\Then{\underbar{Then}}

\define\Dom{\; {\rm Dom} \; }
\define\mod{\; {\rm mod} \; }

\define\cS{{\cal{S}}}

\def\implies{\Rightarrow}

\def\vare{\varepsilon}

\def\b{\beta}
\def\a{\alpha}
\def\l{\lambda}
\def\th{\theta}
\def\z{\zeta}

\def\rest{\mathord\restriction}

\def\w{\omega}

\def\tp{{\rm tp}}

\def\nacc{{\rm nacc}}

\def\g{\gamma}

\def\sig{\sigma}

\def\del{\delta}
\font\frac=eurb10

\def\k{\kappa}

\def\gu{{\frak U}}
\def\gv{{\frak V}}

\define\qed{\hfill $\square$}

\define\subm{\leq _{\frak K}}
\define\suba{\leq _{\mu,\a}^a}
\define\subb{\leq _{\mu,\a}^b}
\define\subc{\leq _{\mu,\a}^c}
\newcommand{\smallbox}[1]{\leavevmode\thinspace\hbox{\vrule\vtop{\vbox
   {\hrule\kern1pt\hbox{\strut\thinspace{#1}\thinspace}}
   \kern1pt\hrule}\vrule}\thinspace}

\newcount\skewfactor
\def\mathunderaccent#1#2 {\let\theaccent#1\skewfactor#2
\mathpalette\putaccentunder}
\def\putaccentunder#1#2{\oalign{$#1#2$\crcr\hidewidth
\vbox to.2ex{\hbox{$#1\skew\skewfactor\theaccent{}$}\vss}\hidewidth}}

\newcount\skewfactor
\def\mathunderaccent#1#2 {\let\theaccent#1\skewfactor#2
\mathpalette\putaccentunder}
\def\putaccentunder#1#2{\oalign{$#1#2$\crcr\hidewidth
\vbox to.2ex{\hbox{$#1\skew\skewfactor\theaccent{}$}\vss}\hidewidth}}

\renewcommand{\iff}{\ \ \underline{\rm iff}\ \ }

\newcommand{\otp}{{\rm otp}}
\newcommand{\Proof}{\noindent{\sc Proof} \hspace{0.2in}}

\newcommand{\Rang}{{\rm Rang}}
\newcommand{\lesdot}{\mathrel{\mathord{<}\!\!\raise 0.8
pt\hbox{$\scriptstyle\circ$}}}

\newcommand{\GK}{{\frak K}}

\newtheorem{theorem}{Theorem}[subsection]

\newtheorem{definition}[theorem]{Definition}

\newtheorem{claim}[theorem]{Claim}
\newtheorem{conclusion}[theorem]{Conclusion}

\newtheorem{remark}[theorem]{Remark} 
 
\newtheorem{fact}[theorem]{Fact}

\newtheorem{hypothesis}[theorem]{Hypothesis}
\newtheorem{construction}[theorem]{Construction}

\title{ Toward categoricity for classes with no maximal models}
\author{{\bf Saharon Shelah}\thanks{\
The first author's research was partially supported
by `BSF' (USA-Israel); Publication no 635 in the first author's
publication list. Done 12/96-1/97/4-97}\\
Institute of Mathematics, The Hebrew University\\
Department of Mathematics, Rutgers University
\and
{\bf Andr\'es Villaveces}\\
Institute of Mathematics, The Hebrew University\\
Departamento de Matem\'atica, Universidad Nacional\\
Bogot\'a, Colombia\\
}
\date{
printed: \today}

\begin{document}
\maketitle

\begin{abstract}
We provide here the first steps toward a Classification Theory of
Abstract Elementary Classes with no maximal models, plus some mild set
theoretical assumptions, when the class is categorical in some
$\lambda$ greater than its L\"owenheim-Skolem number. We study the
degree to which amalgamation may be recovered, the behaviour of non
$\mu$-splitting types. Most importantly, the existence of saturated
models in a strong enough sense is proved, as a first step toward a
complete solution to the \L o\'s Conjecture for these classes. Further
results are in preparation.
\end{abstract}
\eject

\section*{Annotated Content}

\par \noindent
\S0 Introduction.

\begin{quotation}
[We link the present work to previous articles in the same field, and
provide a large-scale picture in which the results of this paper fit.]
\end{quotation}

\par \noindent
\S1 How much amalgamation is left?
\begin{quotation}
[Following \cite{Sh:88}, we prove that although amalgamation is not
assumed, amalgamation bases are `dense' for our purposes. We also
prove the existence of Universal Extensions over amalgamation bases.]
\end{quotation}

\par \noindent
\S2 Types and Splitting

\begin{quotation}
[We provide here the right notions of type for our context, and study
the behaviour of non $\mu$-splitting.]
\end{quotation}

\par \noindent
\S3 Building the right kind of limits

\begin{quotation}
[We define the classes ${}^+\GK _{\mu,\a}$ and variants, in order to
study in depth different concepts of limit models, and of saturated
models. We prove the existence of a good notion of saturated models.]
\end{quotation}

\par \noindent

\section{Introduction}

We study the categoricity spectrum of abstract elementary
classes, when amalgamation is not assumed a priori, and the only
strong model theoretical assumption is the non existence of maximal
models. This looks to us like quite a natural assumption, and many
classes of models that appear usually in mathematics satisfy it---
while they are not first order, and thus need the expansion of
Classification Theory, to which this work contributes.

Previous work with similar motivation appeared in the papers
\cite{Sh:87a}, \cite{Sh:87b}, \cite{Sh:88}, \cite{MaSh:285},
\cite{KlSh:362}, \cite{Sh:394}, \cite{Sh:472}, \cite{Sh:576}, where
the endeavour of extending Classification Theory to more general
classes of models was started. Of course, some additional assumption
had to be used in each one of those directions. There were {\bf set
theoretical} as well as {\bf model theoretical} assumptions.

Among those set theoretical, the main lines were opened by Makkai and
Shelah in \cite{MaSh:285}, where the existence of compact cardinals
was used, and the Categoricity Spectrum for the corresponding classes
was studied. Then followed the work of Kolman and Shelah
\cite{KlSh:362}, and Shelah \cite{Sh:472}, where the hypothesis was
reduced to that of the existence of measurable cardinals. Along those
lines, the \L o\'s Conjecture is not yet fully solved.

Many of the central difficulties in those papers had to do with
pinning down the right kinds of types (when there is no compactness,
the formula-based definition of types is no longer a good one), and with
proving that the amalgamation property for the class still holds. In
\cite{MaSh:285} compactness was still the central tool, and
the definition of types did not present a problem. The compactness
also eased out in a crucial way the proof of amalgamation as well as
the study of the categoricity spectrum. Of course, the price for the
relative smoothness was high; thence the natural motivation of looking
for results with more modest assumptions: reducing the large cardinal
assumption to the existence of a measurable cardinal. This was worked
out in \cite{KlSh:362} and \cite{Sh:472}. A considerable amount of
work was then needed to pin down a notion of `good' extensions. The
lack of compactness was partially supplied for by the use of
Generalised Ultrapowers (of structures). Their existence uses in a
crucial way the measurability, and was central to the proof of the
Categoricity Theorem there.

Among the model theoretical assumptions, the main references are at
this point \cite{Sh:394}, where the amalgamation property is {\bf the
main} assumption. In this context, an extensive use of various kinds
of Ehrenfeucht-Mostowski models is the central tool for constructing
models in the proofs.

This paper could be thought of as `branching off from \cite{Sh:394}'
(here, the amalgamation property is replaced by the weaker
model theoretical assumption of the non-existence of maximal
models). But this is not a completely accurate description of where
this paper fits in the large picture: our set theoretical assumptions are
definitely stronger that those of \cite{Sh:394}: we use GCH in large
chunks of Card, as well as diamonds and weak diamonds. Nevertheless,
we do not use large cardinals, and in this relative sense, this paper
`improves' \cite{KlSh:362} and \cite{Sh:472}.

We plan to continue along this line of research. The forthcoming paper
\cite{ShVi:648} is the next stage.

We shall make free use of EM-models for abstract elementary classes,
throughout the paper.

We wish to thank Mirna D\v zamonja for her useful comments on some
aspects of this paper.


\section{How much amalgamation is left?}

This first section provides the basic framework for the work ---we
study the extent to which amalgamation may be recovered under our
assumptions, as well as the existence of Universal Extensions. We also
provide the main basic definitions.

\subsection{A word about the hypotheses to be used ---\\ Abstract
Elementary Classes.}

The main {\bf model theoretical} hypothesis at work here is, as
indicated in the title, the non-existence of maximal models in the
class. The main {\bf set theoretical} assumption here is the GCH, or
at least the existence of weak diamonds over the relevant
cardinals.

\bigskip
Additionally, we will assume in many parts of this work that the
classes
\begin{description}
\item[(1)] have a L\"owenheim-Skolem number $LS(\GK )$, and
\item[(2)] that they are {\bf categorical} for some $\l$,
\end{description}
with $\l$ high enough compared with $LS(\GK )$, or at least that the
number of models of cardinality $\l$ in $\GK$ is $< \mu_\GK (\l)$,
modulo isomorphism. $\mu_\GK (\l)$ is often equal to $2^\l$ (in this
case, the assumption is just that the class does not have the maximum
possible number of models in $\l$), but in other cases may be `a bit
less' than $2^\l$. For more details on the relationship between
$\mu_\GK (\l)$ and $2^\l$, the reader is referred to \cite[\S
1]{Sh:576}. There, our $\mu_\GK (\l)$ is called $\mu_{\rm wd}(\l)$;
the definition provided there is much more general than what we need
here; we roughly describe $\mu_\GK (\l)$ as `the covering number for
the weak diamond ideal on $\l$'.

In some portions of the work, certain versions of
$\diamondsuit_{S^{\mu^+}_{{\rm cf}(\mu)}}$ for $\mu \in
[LS(\GK),\lambda)$ are used. The full power of GCH is not really
needed throughout the paper; still it is essential for the proof of
the local character of non $\mu$-splitting of types, a central notion
in this work. Up to some point, the set theoretical assumption GCH
`provides' here what otherwise is missing as model theoretical
assumptions, when we compare our hypotheses to those of \cite{Sh:394}
(specifically, the assumption there that all models in $\GK$ are
amalgamation bases).

\begin{definition}{\bf (Abstract Elementary Classes)}\label{absel}

\begin{description}
\item[(1)] $\GK=(K,\subm)$ is an {\bf abstract elementary class} iff
$\GK$ is a class of
models of some fixed vocabulary $\tau=\tau_{\GK }$ and $\subm$
is a two place relation on $K$, satisfying the following axioms
\begin{description}
\item[Ax 0:] If $M\in \GK$, then all $\tau$-models isomorphic to $M$
are also in $K$. The relation $\subm$ is preserved under isomorphisms,
\item[Ax I:] If $M\subm N$, then $M$ is a submodel of $N$,
\item[Ax II:] $\subm$ is an order on $K$,
\item[Ax III:] The union of
a $\subm$-increasing continuous chain $\vec{M}$ of elements of $\GK$ is
an element of $\GK$,
\item[Ax IV:] The union of
a $\subm$-increasing continuous chain $\vec{M}$ of elements of $\GK$
is the lub of $\vec{M}$ under $\subm$,
\item[Ax V:] If $M_\ell\subm N$ for
$\ell \in \{0,1\}$ and $M_0$ is a submodel of $M_1$, then $M_0\subm M_1$,
\item[Ax VI:]
There is a cardinal $\kappa$ such that for every $M\in\GK $ and $A
\subset |M|$, there is $N\subm M$
such that
$A\subset |N|$ and
$\| N\| \leq \kappa\cdot |A|$. The least
such $\kappa$ is denoted by
${\rm LS}(\GK )$ and called {\bf the L\"owenheim-Skolem number of $\GK $}.
\end{description}

\item[(2)] If $\lambda$ is a cardinal and $\GK $ an abstract
elementary class,
we denote by $\GK _\lambda$ the family of all elements of $\GK $
whose cardinality is $\lambda$. We similarly define $\GK _{<\lambda}$

\item[(3)] Suppose that $\GK $ is an abstract elementary class.
\begin{description}
\item[(a)] $\GK $ is said to have the {\bf joint embedding property} (`JEP')
iff for any $M_1, M_2\in \GK $, there is $N\in \GK $ such that
$M_1,M_2$ are
$\subm$-embeddable into $N$.
\item[(b)] $\GK $ is said to have {\bf amalgamation} iff for all $M_0, M_1, M_2
\in \GK$ and $\subm$-embeddings $g_l:\,M_0\longrightarrow M_l$ for
$l\in\{1,2\}$, there is $N\in \GK $ and $\subm$-embeddings
$f_l:\,M_l\longrightarrow N$ such that $f_1\circ g_1=f_2\circ g_2$.
\end{description}

\item[(4)] For $\GK^1\subset \GK$, let

$(\GK ^1)^{am} = \left \{ M_0\in \GK ^1 \left |
\begin{tabular}{l}
if $M_1$, $M_2\in \GK^1$, $g_1$, $g_2$ are as in (3)(b),\\
then there are $N\in \GK^1$, and $f_1$, $f_2$\\
such that $f_1\circ g_1=f_2\circ g_2$
\end{tabular}\right .\right \}$.

The main point here
is to get the amalgamation {\bf inside} the class.
\end{description}
\end{definition}

\subsection{Density of Amalgamation Bases}

To ease the reading of this paper, we shall (sometimes redundantly)
endeavour to spell out the hypotheses used, at the beginning of each
section.

\begin{hypothesis}\label{nomax}
$\GK$ is an abstract elementary class with no maximal model in $\GK
_{<\l}$, categorical in $\l$.
\end{hypothesis}

The content and the proof of Theorem \ref{density} below are basic in
subsequent work. They are akin to those in \cite[Theorem 1.3]{Sh:87a} and
\cite{Sh:88}. Still, for the sake of completeness, we provide the
argument. See more in \cite{Sh:576}, for a study of
Weak Diamond principles and their relation to model theoretical
properties of Abstract Elementary Classes.

\bigskip

Before looking at the next central questions in this context, namely
the density of amalgamation bases and the existence of universal
extensions over every model in $\GK_\k$, we need some additional
results.

\begin{fact}\label{jep}
Suppose that there are no maximal models in $\GK _{\leq \l}$. Fix
cardinals $\chi$, $\mu$ such that $LS(\GK ) \leq \mu \leq \chi$. Then
\begin{description}
\item[1)] if $M \in \GK_{< \chi}$ and $\|M\| + LS(\GK ) \leq
\mu \leq \chi$, then there is $N$ such that
$M \subm N \in \GK_\mu$.
\item[2)] If $\GK$ is categorical in $\chi$ (in particular, if $\chi
=\l$), then $\GK_{\leq \chi}$ has the JEP (joint embedding property).
\end{description}
\end{fact}

\Proof 1) is easy, from a repeated use of the L\"owenheim-Skolem
theorem (in this context), and the non existence of maximal models:
axiom III guarantees that unions of $\subm$-increasing continuous
chains of elements of $\GK$ are in $\GK$.

\noindent
2) is also easy to see, by embedding the models into extensions of
size $\chi$.
\qed$_{\ref{jep}}$
\bigskip

The main tool to construct models which have useful homogeneity
properties is in this context the use of generalised
Ehrenfeucht-Mostowski models. These were developed for the context
of abstract elementary classes by Saharon Shelah in \cite{Sh:88}. The
following fact asserts that they exist in this context.

\begin{fact}\label{emmodels}
For every linear order $I$, there is $\Phi$ such that $EM(I,\Phi)$
is an EM model (so, for instance, if $EM(I,\Phi)\in \GK$ and $J\subset
I$, then $EM(J,\Phi)\subm EM(I,\Phi)$).
\end{fact}

\Proof Since there are no maximal models in $\GK$, there are models in
$\GK_\mu$, where $\mu= |EM(I,\Phi)| = |I| + |\tau| + \beth_{(2^{LS(\GK
)})^+}$, by \cite[1.7]{Sh:88} (where $\tau$ is the size of the
vocabulary). Now the construction of the EM models can be carried in a
way similar to how it was done in \cite[VII,\S5]{Sh:c}.

\qed$_{\ref{emmodels}}$

\begin{theorem}\label{density}{\bf (Density of Amalgamation Bases)}
If $LS(\GK) < \k \leq \l$ (remember: $\l$ is the categoricity cardinal
of the hypotheses), $\exists \theta (2^\theta = 2^{<\k} <
2^\k)$, then for every $M \in \GK_{<\k}$, there is $N$ with $M\subm N
\in \GK^{am}_{<\k}$.
\end{theorem}

\Proof Suppose $M$ is a counterexample to this. The idea is to build a binary
tree of models on top of $M$, in such a way that the two immediate
successors of every node act as counterexamples to amalgamation over
$M$, and then use the weak diamond at $\k$ (whose existence is
guaranteed by $2^\th<2^\k$! ---for more on generalised weak diamonds,
see the Appendix to the forthcoming `Proper Forcing' book by the first
author \cite{Sh:f}) to get a contradiction. So, we choose by induction
on $\a<\k$ models $M_\eta$, for $\eta\in {}^\a 2$, such that

\begin{description}
\item [(a)]  $M_{<>} = M$
\item [(b)]  $M_\eta \in \GK_{< \kappa}$
\item [(c)]  $\alpha$ limit $\wedge$ $\eta \in {}^\alpha 2 \Rightarrow M_\eta
=  \bigcup\limits_{\beta < \alpha} M_{\eta \restriction \beta}$
\item [(d)]  $\beta < {\rm lg}(\eta) \Rightarrow M_{\eta \restriction
\beta} \subm M_\eta$
\item [(e)]  $M_{\eta \peg \langle 0 \rangle},
M_{\eta \peg \langle 1 \rangle}$ cannot be amalgamated over $M_\eta$;
i.e. there is no $N \in \GK_{<\kappa}$ and $\subm$-embedding $f_\ell:
M_{\eta \peg \langle \ell \rangle} \rightarrow N$ such that
$f_0 \restriction M_\eta = f_1 \restriction M_\eta$
(so $M_\eta \not= M_{\eta \peg \langle \ell \rangle}$).

\end{description}
\medskip

\noindent
For each $\eta \in {}^\kappa 2,M_\eta = \bigcup\limits_{\alpha <
\kappa} M_{\eta \restriction \alpha} \in \GK_\kappa$ hence by
Fact~\ref{jep}(1) there is $N_\eta
\in \GK_\lambda$ with $M_\eta \subm N_\eta$.  By the categoricity in $\l$,
there exists an isomorphism $h_\eta:N_\eta {}\stackrel{\rm
onto}{\longrightarrow} N^* := EM(\lambda,\Phi)$.  But then
$h_\eta(M_\eta) \subm N^*$, hence $M_\eta$ is $\subm$-embedded into
$EM(\alpha_\eta,\Phi)$ for some $\alpha_\eta < \kappa^+$. 

Let $<^*$ be a linear order on $\k$ isomorphic to
$({}^{\w>}\k,<_{lex})$, so that each $\a<\k^+$ can be embedded into
it. Then $EM(\a_\eta,\Phi)$ is $\subm$-embeddable into
$N^*=EM((\k,<^*),\Phi)$. So, there is a $\subm$-embedding
$h_\eta^*:M_\eta \longrightarrow N^*$.

Now use the weak diamond: since there exists
$\theta$ such that $2^\theta = 2^{<\k} < 2^\k$, the weak diamond for
$\k$ holds, and thus there are distinct $\eta_1$, $\eta_2\in {}^\kappa
2$ and there is $\a<\k$ such that $h^*_{\eta _1}\rest M_{\eta _1\rest
\a} = h^*_{\eta _2}\rest M_{\eta _2\rest \a}$, and
$\eta_1(\a)\not=\eta_2(\a)$. But both $M_{\eta _1\rest
{\a+1}}$ and $M_{\eta _2\rest {\a+1}}$ embed into
$EM(\k^+,\Phi)$. This contradicts that $M_{\eta _1\rest \a}$ is not an
amalgamation base!

\qed$_{\ref{density}}$
\bigskip

So, we have density of amalgamation bases in the case mentioned above
(there exists $\theta$ such that $2^\theta = 2^{<\k} < 2^\k$), but it
should be made clear that the use of the weak diamond (or, a fortiori,
of GCH in $[LS(\GK),\l)$), was crucial here.

\bigskip

\subsection{Universal Extensions}

At this point, we begin to include the following assumption:

\bigskip
\noindent
{\bf GCH Hypothesis}: $2^\mu=\mu^+$, for all
$\mu\in[LS(\GK),\l)$. Although we stated at the outset this
assumption, we repeat it now. Up to now, the weak diamond was enough.
Nevertheless, it is worth stressing that
our aim is to obtain as much stability as possible for our new
contexts, and at the same time trying to use as little as possible
set-theoretical assumptions. GCH does not seem too unreasonable from
this point of view.

\bigskip

The following theorem is crucial in the study of the right kind of
types in our context, and is a natural step in allowing us to build
models with enough saturation. So far, we have not defined the types,
and thus we concentrate on universality. It is worth noting that the
existence of universal extensions here is obtained {\bf for
amalgamation bases.}

\begin{theorem}\label{universal}{\bf (Existence of Universal
Extensions)}
Suppose that $\mu \in [LS(\GK),\lambda)$ and $M_0 \in
\GK^{am}_\mu$. Then there is $M_1$ such that $M_0 \subm M_1 \in
\GK^{am}_\mu, M_1$ is universal over $M_0$ (i.e. $M_0 \subm M_2 \in
\GK_\mu \Rightarrow M_2$ is $\subm$-embeddable into $M_1$ over $M_0$).
\end{theorem}

\Proof Let $I$ be a linear order of cardinality $\mu^+$ such that $I \times
(\a^+ + 1) \approx I$, for every $\a < \mu^+$, and pick $M_0 \in
\GK^{am}_\mu$. We first move to the case of EM models, and prove the
following fact.

\begin{claim}\label{useofgch}
There is a $\subm$-embedding $f:M_0\longrightarrow EM(I,\Phi)$
such that for every $M_1$ with $M_0\subm M_1 \in \GK _\mu$
there is a $\subm$-embedding $g:M_1 \longrightarrow EM(I,\Phi)$ extending $f$.
\end{claim}

\Proof We begin by listing (note the strong use of $2^\mu=\mu^+$ here!)
all the possible embeddings from $M_0$ into $EM(I,\Phi)$ as $\langle f_i |
i<\mu^+\rangle$. For every $f_i$ let now $M_{1,i}\in \GK_\mu$ be a
counterexample to the property we are looking for; namely, $M_0\subm
M_{1,i}$ and
$f_i$ does not `lift' to an embedding from $M_{1,i}$ to
$EM(I,\Phi)$. Since $M_0\in \GK^{am}_\mu$, we can find $\langle
M_{2,i} | i\leq \mu^+\rangle$, $\subm$-increasing, continuous, such
that $M_{2,0}=M_0$, and $i<\mu^+ \implies M_{1,i}$ is embeddable into
$M_{2,i+1}$ over $M_0$. Now, by categoricity in $\l$, we know that the
limit $M_{2,\mu ^+}$ is embeddable into $EM(\l,\Phi)$, and thus into
some $EM(\a^*,\Phi)$, $\a^* < \mu ^{++}$, and hence into $EM(I,\Phi)$,
say by $g$. But then $g\rest M_0$ must be $f_{i(*)}$, for some
$i(*)<\mu^+$. Contradiction. \qed$_{\ref{useofgch}}$

\begin{figure}[ht]\hspace*{35.00mm}
\psfig{figure=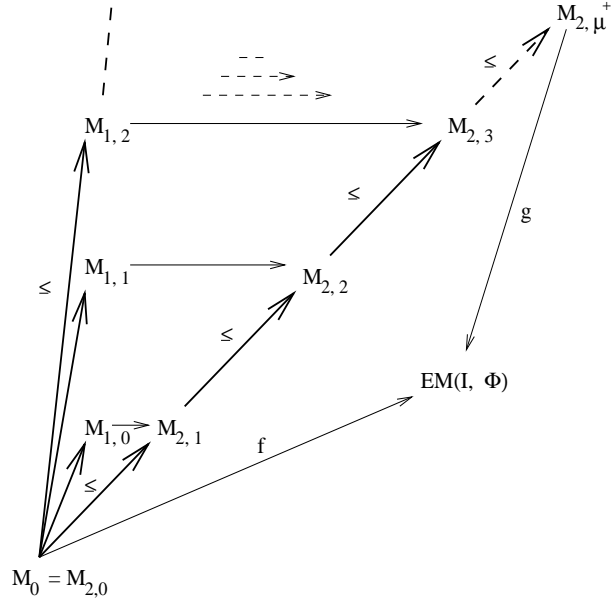,height=8cm,width=8cm}
\caption{ Lifting $f$, in Claim \ref{useofgch} } 
\end{figure}

\bigskip
\noindent
Fix now some $f$ as in the claim, and let $I_0\subset I$, $|I_0|= \mu$
be such that Rang($f$)$\subset EM(I_0,\Phi)$. We could have chosen $I$
from the beginning as being decomposable as
$$I=\bigcup\limits_{\zeta<\mu^+}I_\zeta,$$
for $\langle
I_\zeta\rangle_{\zeta<\mu^+}$ increasing, $|I_\z|=\mu$, and $I_{\z
+1}$ universal (inside $I$) over $I_\z$ (a similar construction is
also used in \cite{Sh:220}, and may be obtained by taking e.g.
$$I=\{
\eta\in {}^\omega (\mu^+) | \eta \mbox { eventually 0 but not
constantly } 0\},$$
ordered lexicographically, and $I_\z = \{ \eta \in I |
\Rang (\eta) \subset \mu \times (1+\z )\}$).

\bigskip
\noindent
Let now $M_1^*$ correspond via an isomorphism $f^+ \supset f$,
$$f^+:M_1^*\stackrel{{\rm onto}}{\longrightarrow}
EM(I_1,\Phi),$$
to $EM(I_1,\Phi)$. We claim that $M_1^*$ is universal over $M_0$: Let
$M_0\subm M_2\in \GK_\mu$. Pick the corresponding embedding
$f_2\supset f$, $f_2:M_2 \longrightarrow EM(I,\Phi)$. As before, let
$I'\subset I$, $|I'|\leq \mu$ be such that $\Rang (f_2)\subset
EM(I',\Phi)$. Thus, for some automorphism $h$ of $I'$, such that
$h\rest I_0 =$ identity, $h``(I')\subset I_1$. Then $h$ induces an
automorphism $\hat h$ of $EM(I,\Phi)$, and we have
$$M_2 \stackrel{f_2}{\longrightarrow} EM(I',\Phi) \stackrel{\hat
h}{\longrightarrow} EM(I_1,\Phi) \stackrel{(f^+)^{-1}}{\longrightarrow}
M_1^*,$$ and $\hat h \rest \Rang (f) =$ id, $f_2\supset f$. So,
$(f^+)^{-1}\circ \hat h \circ f_2$ is an isomorphism from $M_2$ into
$M^*_1$, and its restriction to $M_0$ is the identity.
\qed$_{\ref{universal}}$

\bigskip

Thus, we have universal models in the right cardinals over
amalgamation bases. The following definition should be regarded as a
first step toward the $(\mu,\nu)$-limits and our version of
saturation.

\begin{definition}
Let $M_0 <^2_\mu M_1$ mean (for $\mu \in [LS(\GK),\lambda)$) that
$M_0,M_1\in \GK_\mu$, $M_0\subm M_1$ and $M_1$ is universal over $M_0$.
\end{definition}

\bigskip

\begin{definition}\label{<3}
Let $\delta$ be a limit ordinal, $\del \leq \mu^+$,
$\mu \in [LS(\GK),\lambda)$. Then,
$$M <^3_{\mu,\delta} N$$ $\iff$ there is a $\subm$-increasing continuous
sequence $\vec{M} = \langle M_i|i \le \delta \rangle$ such that $M_0 =
M$ and $i<\del \implies M_i \in \GK^{am}_\mu$ and $M_i <^2_\mu
M_{i+1}$, and $N = \bigcup_{i < \delta} M_i$ (so $\|N\| = \mu +
|\delta|$).
\end{definition}

\noindent
(In this case, we use $M_\delta$ for $\bigcup_{i < \delta} M_i$ and
call $\langle M_i: i \le \delta \rangle$ a {\bf witness for $M
<^3_{\mu,\delta} N$}.)
\bigskip

\begin{remark}
\end{remark}
\noindent
In previous uses of these extensions, the amalgamation property was
assumed to hold in the class --- here we must stress the fact that
{\bf by decree} all the levels up from $M$ to $N$ {\it are
amalgamation bases}.

\bigskip
Among the basic properties of $<^3_{\mu,\delta}$, we have that

\begin{fact}\label{3andlimits}{\bf ($<^3_{\mu,\delta}$ and limits)}

\begin{description}
\item[1)] If $M_\ell <^3_{\mu,\delta} N_\ell$ for
$\ell =1,2$ and $h$ is a $\subm$-embedding of $M_1$ into $M_2$
\underbar{then} we can extend $h$ to an isomorphism $h^*$ from 
$N_1$ onto $N_2$.
\item[2)] Moreover, if $\langle M_{\ell,i}:i \le \delta \rangle$ witnesses
$M_\ell <^3_{\mu,\delta} N_\ell$ we can demand that $h^*$ map $M_{1,2i}$ into
$M_{2,2i}$ and $(h^*)^{-1}$ map $M_{2,2i+1}$ into $M_{1,2i+1}$.
\end{description}
\end{fact}
\Proof By induction on $\del$.
\qed$_{\ref{3andlimits}}$

\bigskip
Another easy fact about $<^3_{\mu,\delta}$ is:

\begin{fact}
\begin{description}
\item[1)] $M <^3_{\mu,\delta}N$ \iff $M <^3_{\mu,\cf(\delta)}N$.
\item[2)] If $M <^3_{\mu,\delta}N$ and $u\subseteq\del=\sup(u)$,
$\del_1=\otp(u)$, then $M <^3_{\mu,\delta_1}N$.
\end{description}
\end{fact}

\bigskip

\begin{hypothesis}  $\diamondsuit_{S^{\mu^+}_{{\rm cf}(\mu)}}$
for $\mu \in [LS(\GK),\lambda)$.
\end{hypothesis}

The use of amalgamation wherever possible, together with the existence
of universal models over amalgamation bases, are the two basic tools
of construction of saturated enough models.  The following fact is
important from that point of view.

\bigskip
\begin{fact}\label{amalg3}
If $M <^3_{\mu,\delta}N$ and $\del<\mu^+$, then $N\in\GK^{am}_\mu$.
\end{fact}

\Proof Easy by $\diamondsuit_{S^{\mu^+}_{{\rm cf}(\mu)}}$ (proof
similar to that of \ref{density}).

\qed$_{\ref{amalg3}}$

\begin{fact}\label{3andlimits2}
For every limit $\del<\mu^+$ we have
\begin{description}
\item[1)] If $M\in \GK ^{am}_\mu$, then for some $N$ we have
$M<^3_{\mu,\del}N$,
\item[2)] If $M<^3_{\mu,\del}N$, then $N\in \GK ^{am}_\mu$.
\end{description}
\end{fact}
\Proof By induction on $\del$. Suppose that this is true for all limit
ordinals $<\del$. If $\del$ is not a limit of limits, let $\del_0$ be
the highest limit below $\del$. We have by induction hypothesis $M_{\del_0}$
with $M<^3_{\mu,\del_0}M_{\del_0}$, as witnessed by some sequence
$\langle M_i | i\leq \del_0\rangle$. Just taking a universal extension
$M_\del$ of $M_{\del_0}$
over $M$ does the trick: by Fact~\ref{amalg3}, $M_{\del_0}$ is itself an
amalgamation base, and thus the sequence $\langle M_i | i\leq
\del_0\rangle\peg \langle M_\del\rangle$ witnesses that
$M<^3_{\mu,\del}M_\del$.

Now, if $\del$ is a limit of limits, we construct a $\subm$-increasing
continuous sequence $\langle M_i | i\leq
\del\rangle$ of models {\bf from $\GK_\mu$} such that $M_{i+1}\in \GK
^{am}$, $M_{i+2}$ universal over $M_{i+1}$ over $M = M_0$. We just take
unions at limits, and universal extensions which are amalgamation
bases at successors. This is close to what we need, but of course we
still need to get that the $M_i$'s are amalgamation bases {\bf all the
way through} $\del$. But this is taken care of by Fact \ref{amalg3}.


\qed$_{\ref{3andlimits2}}$

\bigskip
\noindent
The following definition is crucial in the study of saturatedness in
this class, and will play a central role from now on.

\begin{definition}
We say that $N$ is a {\bf $(\sig,\delta)$-limit} if $M
<^3_{\sig,\delta} N$ for some $M$ and $\sig\leq\del$ is regular.
\end{definition}

The proofs of the previous facts essentially depended on constructions
by induction on $\delta$, as well as the use of set theoretical
hypotheses. These hypotheses are not too strong (from our point of
view), especially when one compares them to those that were used in
the past by Makkai and Shelah in \cite{MaSh:285} (compact cardinals)
and by Kolman and Shelah in \cite{KlSh:362} and Shelah in
\cite{Sh:472} (measurable cardinals). Of course, the `price to pay' is
that many notions `natural' in those contexts (such as definitions of
types as sets of formulas in the presence of a strongly compact
cardinal in \cite{MaSh:285} or the existence of many ultrapower
operations in \cite{KlSh:362} and \cite{Sh:472}) are no longer
`natural' here, and require new ways of dealing with the categoricity
problem.

\bigskip

\section{Types and Splitting}

\subsection{What kind of types are good here?}

We start by giving a definition of types for this context. It must be
stressed that here, types are only defined over {\bf models which are
amalgamation bases}, so as to avoid confusion later. The definition of
types here is
essentially the same from \cite[Ch.II]{Sh:300} and \cite[\S0]{Sh:576}
and \cite{Sh:394}. There is, though, a difference: in the presence of
monster models (like in \cite{Sh:394}), it is natural to construe all the
automorphisms relevant to the definition of types as automorphisms of
the monster. Here, in its absence, we must do with embeddings into an
amalgam. Our hypotheses about amalgamation clear away the problem
here. Still, the diagram chasing involved might be slightly more
entangled than within monster models.

\begin{definition}
\end{definition}
\begin{description}
\item[1)] We define the {\bf type}
$$\tp(\bar a,M,N)$$
(when
$M \subm N$, $\bar a \subset N$, $M,N\in \GK^{am}_\mu$) 
as $(\bar a,M,N)/E$ where $E$ is the following equivalence relation:
$(\bar a^1,M^1,N^1) E (\bar a^2,M^2,N^2)$ \iff $M^\ell
\subm N^\ell,\bar a^\ell \in {}^\alpha(N^\ell)$ (for some $\alpha$),
$M^1 = M^2$ and there is $N \in \GK$ satisfying $M^1 = M^2 \subm
N$ and $f^\ell:N^\ell \stackrel{\subm}{\longrightarrow} N$ over $M^\ell$
(i.e. $f \restriction M^\ell$ is the identity) and $f^1(\bar a^1) =
f^2(\bar a^2)$.

More generally, for $N\in \GK_\mu$ (not necessarily an
amalgamation base) and $\bar a \subset N$, we define $\tp(\bar a,M,N)$
as $\tp(\bar a,M,N')$, with $M\subm N'\subm N$, $\bar a \subset N'$,
and $N'$ is an amalgamation base.
\item[2)] We say that $N$ is $\kappa$-{\bf saturated}
(when $\kappa > LS(\GK)$) if $M \subm N, |M| < \kappa$ and $p \in
{\cal S}^{< \omega}(M)$ (see below) imply that $p$ is realized in $M$,
i.e. for some $\bar b \subset N$, $p = \tp(\bar b,M,N)$.
\item[3)] ${\cal S}^\alpha(M) := \{ \tp(\bar a,M,N) |\bar a \in
{}^\alpha N, M \subm N\}$.
\item[4)] When $M \subm N$ and $p
\in {\cal S}^\a(N)$, we denote by $p \restriction M$ the {\bf
restriction} to $M$ of $p$ given by  $\tp(\bar a,M,N_1)$,
where $N \subm N_1$, $p =
\tp(\bar a,N,N_1)$ and $p \le q$.
\item[5)] ${\cal S}(M) = {\cal S}^1(M)$ (we could just as well use
${\cal S}^{< \omega}(M)$).
\end{description}

\bigskip
\begin{remark}
\end{remark}
We define types on $M$ in $N$ under the condition that $M$ be an
amalgamation base and there be some amalgamation base $N'\supset \bar
a$ in between $M$ and $N$. Under these conditions, we may prove that
$E$ is an equivalence relation. Otherwise, the diagram chasing for the
transitivity of $E$, which we leave to the reader, would not go through.

\bigskip

The following fact is basic, and is used throughout the paper.

\begin{fact}\label{stabbelow}{\bf (Stability below $\l$)}
Let $\mu<\l$. Since $\GK$ is categorical in $\l$, for every $N\in \GK
_\mu$, $|\cS (N)|\leq \mu$.
\end{fact}

\Proof Use $EM(\l,\Phi)$, just as in \cite[Theorem 3.9]{KlSh:362}, where
Kolman and Shelah prove the existence of weakly universal models over
any $N\in \GK_\mu$.
\qed$_{\ref{stabbelow}}$

\begin{definition}
\begin{description}
\item[(1)] {\bf ($\mu$-splitting)} $p \in {\cal S}(M)\,\,\mu$-splits over $N \subm M$ \iff
$|N| \le \mu$, and there are $N_1,N_2,h$ such that: $N_1,N_2\in \GK_\mu$,
$h$ an elementary mapping from $N_1$ onto $N_2$ over $N$ such that
the types $p \restriction N_1$ and $h(p\restriction N_1)$ are contradictory
and $N \subm N_\ell \subm M$.
\item[(2)] We say that the type $q\in {\cal S}(N)$ is a
{\bf stationarisation} of $p\in {\cal S}(M)$, $\|M\|=\mu$, $M\subset
N$, \iff for some $M^-\subm M$, $q$ does not $\mu$-split over $M^-$.
\end{description}
\end{definition}

The next theorem marks the real beginning of the new ideas in this
paper. It uses GCH in a rather strong way, and sheds light on the
local character of non-$\mu$-splitting.

\subsection{The Splitting}

\bigskip
\begin{theorem}\label{splitting}
Assume that
\begin{description}
\item[(a)] $\langle M_i | i\leq\sig\rangle$ is $\subm$-increasing and
continuous,
\item[(b)] for all $i\leq\sig$, $M_i\in\GK_\mu$ and $M_i$ is an
amalgamation base in $\GK_\mu$,
\item[(c)] each $M_{i+1}$ is universal over $M_i$,
\item[(d)] $\cf (\sig)=\sig\leq\mu^+\leq \l$, and $p\in\cS(M_\sig)$. (Since
$M_\sig$ is an amalgamation base, `types' are well-defined in this
context.)
\end{description}
\Then, for some $i<\sig$, $p$ does not $\mu$-split over $M_i$.
\end{theorem}

\begin{remark}
\end{remark}
We do not just have
$$\bigwedge_{j\in (i,\sigma)}[p\rest M_j \mbox { does not
split over } M_i].$$

\Proof Assume that the conclusion fails. We shall choose $\langle M_i
| i\leq\sig\rangle$ and $p$ contradicting
the statement, fitting into one of the following possibilities.

\begin{description}
\item[(a)] $j<\sig \implies p\rest M_j$ does not $\mu$-split over $M_0$.
\item[(b)] else (a) is impossible, and $p\rest M_{2i+1}$
$\mu$-splits over $M_{2i}$ and $p\rest M_{2i+2}$ does not $\mu$-split
over $M_{2i+1}$.
\item[(c)] else (a) and (b) are both impossible, and
$\sig=\mu$ (so $\mu$ is regular), and $i<\sig \implies p\rest M_{i+1}$
$\mu$-splits over $M_i$.
\end{description}

Without loss of generality, $M_{i+1}$ is $(\mu,\w)$-limit over $M_i$,
as there is such an $M_{i+1}'$, $M_i\subm M_{i+1}'\subm M_{i+1}$.

\begin{claim}\label{split}
One of (a), (b) and (c) is always possible.
\end{claim}

\Proof Assume that both (a) and (b) are impossible. Given
$\vec{M}=\langle M_i | i\leq\sig\rangle$ and $p\in\cS(M_\sig)$, we
will use the fact that both (a) and (b) are impossible (for any $\vec{M}$)
in order to produce some $\vec{M}'$ satisfying (c). As for any
$j<\sig$, possibility (a) fails for $\langle M_{j+i}|i\leq
\sig\rangle$, we have that necessarily

\begin{description}
\item[$(*)$] for every $j<\sig$, there is some $\zeta_j\in (j,\sig)$
such that $p\rest
M_{\zeta_j}$ $\mu$-splits over $M_j$.
\end{description}
\noindent
Even more so, by renaming, we can require
\begin{description}
\item[$(*)'$] $p\rest M_{i+1}$ $\mu$-splits over $M_i$. [We are
close here to a $(c)$-style sequence. What is still missing is the
appropriate length.]
\end{description}

We can find $\langle M_{i,j} | j\leq\mu\rangle$ $\subm$-increasing
continuous, $M_{i,j+1}$ $(\mu,\w)$-limit over $M_{i,j}$, $M_{i,j}$ an
amalgamation base [we freely use \ref{3andlimits2}], for each
$j\leq\mu$, with $M_{i,0} = M_i$, $M_{i,\mu}\subm M_{i+1}$. Now we
ask, for each $i$,
\begin{description}
\item[$\otimes_i$] Does $p\rest M_{i,\mu}$ $\mu$-split over $M_{i,j}$
for every $j\leq \mu$?
\end{description}
If for some $i$, the answer is `yes', then we can repeat the
procedure above
(applied now to $\langle M_{i,j} | j\leq\mu\rangle$ and $p\rest
M_{i,\mu}$). So we get that $(*)'$ holds, i.e. possibility (c)
holds for $\langle M_{i,j}|j\leq \mu\rangle$. If, on the other hand,
for every $i\leq \sig$, the answer to
$\otimes_i$ is no, then for some $j_i<\mu$, $p\rest M_{i,\mu}$ does
not $\mu$-split over $M_{i,j_i}$.  Consider the sequence
$$\langle
M_{0,j_0},M_{0,\mu},M_{1,j_1},M_{1,\mu},\dots
,M_{\sig,j_\sig},M_{\sig,\mu}\rangle.$$
This sequence and $p$ clearly witness case $(b)$.
\qed$_{\ref{split}}$

\bigskip
We now come back to the proof of Theorem \ref{splitting}, and look at
the three possible cases from the last claim.

\bigskip
\noindent
{\bf Proof in Case (c):} Under our hypotheses, we have that $2^{<\mu} =
\mu$. Let $p$, $\langle M_i | i\leq
\mu\rangle$ be as in case (c). Choose by induction on $i\leq \mu$ models
$N_i$ and sequences $\langle g^i_\eta | \eta \in {}^i 2\rangle$ such that

\begin{description}
\item[$(\a)$] $N_i\in \GK_\mu$, $\langle N_i | i\leq\mu\rangle$ is
$\subm$-increasing continuous,
\item[$(\b)$] $N_{i+1}$ is a $(\mu,\omega)$-limit over $N_i$,
\item[$(\g)$] $N_0 = M_0$,
\item[$(\del)$] $g^i_\eta$ is an isomorphism from $M_{\omega i}$ onto $N_i$,
\item[$(\varepsilon)$] $g^i_{\eta\rest j}\subset g^i_\eta$, for $j<i$,
\item[$(\zeta)$]$g^{i+1}_{\eta\peg <0>}(p\rest M_{\omega (i+1)})\not=
g^{i+1}_{\eta\peg <1>}(p\rest M_{\omega (i+1)})$.
\end{description}

The clause $(\zeta)$ is possible because $\vec{M}$, $p$ witness the
case (c).
Having obtained these sequences and isomorphisms, we have that $N_\mu$
is a $(\mu,\mu)$-limit. For $\eta\in {}^\mu 2$, $g_\eta :M_\mu
\stackrel{\approx}{\longrightarrow} N_\mu$, and the $g_\eta (p)$, for
$\eta\in {}^\mu 2$ (in $\cS (N_\mu)$) are pairwise distinct.

So, $N_\mu\in\GK_\mu$ is an amalgamation base, and $|\cS
(N_\mu)|>\mu$. This contradicts the basic fact $\ref{stabbelow}$, and
ends the proof when dealing with possibility (c).

\bigskip
\noindent
{\bf Proof in Case (a) or (b):} Choose $\vec{C}=\vec{C}^\sig = \langle
C^\sig_\a | \a\in S^{\mu^+}_\sig\rangle$, where
\begin{description}
\item[(a)] $S^{\mu^+}_\sig$ denotes the set of ordinals $<\mu^+$ of
cofinality $\sig$, and
\item[(b)] for every $\a$, $C^\sig_\a\subset \a$ is a club, $\otp
(C^\sig_\a) = \sig$, and for every club $C$ of $\mu^+$, the set
$$\Big\{ \del\in S^{\mu^+}_\sig \Big| \del = \sup [C\cap \nacc
(C_\del)]\Big\}$$ is stationary, where $\nacc(X)$ is the set of
nonaccumulation points of $X$. This is possible by \cite[III]{Sh:g}.
\end{description}

\bigskip
We start with $p$ and $\langle M_i | i\leq
\sig\rangle$ as there and choose (by induction on $\a<\mu^+$) $N_\a\in
\GK^{am}_\mu$ such that
\begin{description}
\item[(i)] $\langle N_\a |\a<\mu^+\rangle$ is $\subm$-increasing continuous,
\item[(ii)] $N_{\a+1}$ is $(\mu,\w)$-limit over $N_\a$,
\item[(iii)] when $\cf (\a) =\sig$, then we list $C^\sig_\a$ (our
originally chosen club in $\a$ of order type $\sig$) increasingly as
$$C^\sig_\a = \{ \b_{\sig,\a,\zeta} | \zeta<\sig\}.$$
Additionally, we let
$\b_{\sig,\a,\sig} =\a$ and also let $\langle M_i | i\leq
\sig\rangle$ and $\langle N_{\b_{\sig,\a,\zeta}} | \zeta\leq \sig
\rangle$ be isomorphic via $g_\a:M_\sig
\stackrel{\approx}{\longrightarrow} N_\a$ (so that $g_\a(M_\zeta)
= N_{\b_{\sig,\a,\zeta}}$). Let $a_\a \in N_{\a+1}$ realise $g_\a(p)$.
\end{description}

\noindent
So, we have $\langle N_\a|\a<\mu^+\rangle$. Let $N = \bigcup_{\a<\mu^+}N_\a \in
\GK_{\mu^+}$.

\bigskip
Clearly, $N$ $\subm$-embeds into $EM(\l,\Phi)$. Even
more, we can use $\Phi'=\Psi\circ\Phi$ such that $EM(\mu^+,\Phi')$ is
universal in $\GK_{\mu^+}$, and has as many automorphisms as we will
need. For more details on the theory of EM models for abstract
elementary classes, see \cite[I, \S 4]{Sh:394}.

\noindent
So we have a $\subm$-embedding $h:N\longrightarrow
EM(\mu^+,\Phi')$. For $\a\in S$ ($S := S^{\mu^+}_\sig$; $S=\Dom
(\vec{C}^\sig)$), let
$$h(a_\a)=\tau_\a(\xi^\a_1,\dots,\xi^\a_{n(\a)}),$$ with
$\xi^\a_1<\dots<\xi^\a_{m(\a)}<\a\leq\xi^\a_{m(\a)+1},\dots,\xi^\a_{n(\a)}<\mu^+$,
and let 
$$E=\Big\{ \g<\mu^+ \Big| \forall c\in N[c\in N_\g \leftrightarrow h(c)\in
EM(\g,\Phi ')]\mbox { and }\g\mbox { a limit ordinal}\Big\} .$$ Clearly, $E$
is a club.

\bigskip
\noindent
We now focus on case (a): for some stationary $S^*\subset S$, $\a\in
S^*\implies C^\sig_\a \subset E$, $\tau_\a = \tau^*$, $n(\a)=n^*$,
$m(\a)=m^*$, $\xi^\a_1 = \xi^*_1,\dots,\xi^\a_{m(\a)}=\xi^*_{m^*}$,
$\b_{\a,\sig,0} =
\b_{*,0}$. Let $\a'<\a''$ be in $S^*$. We then have that

\bigskip
$\quad \left \{
\begin{tabular}{l}
(i) $\tp(a_{\a''},N_{\a'},N_{\a''+1})$ does not $\mu$-split over $N_{\b_*,0}$\\
(ii) $\tp(a_{\a'},N_{\a'},N_{\a''+1})$ $\mu$-splits over $N_{\b_*,0}$.

\end{tabular}\right .$

\medskip
\noindent
For (i), we use the choice of $p$ and $\langle M_i | i\leq \sig
\rangle$ as in case (a) and choose $j<\sig$ such that
$g_{\a''}(M_j)\supset N_{\a'}$; since $a_{\a''}$ realises
$g_{\a''}(p)$, we get that $\tp(a_{\a''},N_{\a'},N_{\a''+1})$ does not
$\mu$-split over $N_{\b_*,0}$. To see (ii), we just use our original
assumption about the splitting of $p$, and `translate' it via $g_{\a'}$.

\bigskip
\noindent
So, the two types must be different, and thence
\begin{description}
\item $\tp(h(a_{\a''}),EM(\a',\Phi'),M_{\a''+1})\not=
\tp(h(a_{\a'}),EM(\a',\Phi'),M_{\a''+1})$,
\end{description}
but on the other hand, it is easily seen that $h(a_{\a''})$ and
$h(a_{\a'})$ realise the same type ---$\Phi'$ could have been chosen at
the outset so
that there is an automorphism $k$ of $EM(\mu^+,\Phi')$ with
$k\rest EM(\a',\Phi')=$ identity and
$k(h(a_{\a''}))=h(a_{\a'})$.

\bigskip
\noindent
We now switch to case (b): Let $\chi$ be large enough, and let
$\langle {\frak B}_\a |
\a<\mu^+\rangle$ be a $\subm$-increasing continuous sequence of
elementary submodels of $({\cal H}(\chi),\in,<^*_\chi)$, each ${\frak
B}_\a$ of size $\mu$, such that
$\Phi$, $EM(\l^+,\Phi)$, $h$, $\langle M_\a | \a<\mu^+\rangle$ and
$\langle a_\a | \a\in S\rangle$ all belong to ${\frak B} _0$, $\langle
{\frak B}_\a | \a\leq\g\rangle \in {\frak B}_{\g +1}$, and ${\frak
B}_\g\cap \mu^+$ is an ordinal. Let
$$E^*=\{ \g | {\frak B}_\g \cap \mu^+ = \g\}.$$
$E^*$ is a club of $\mu^+$. Also, by the choice of $\vec{C}$, there is
$\a$ such that $C^\sig_\a \subset E^*$. Now find $\zeta<\sig$
such that $\xi^\a_1,\dots,\xi^\a_{m(\a)}<\b_{\sig,\a,\zeta}<\a$, and $p\rest
M_{\zeta+1}$ does not $\mu$-split over $M_\zeta$.

Let now $\varphi$ be a formula in the language of set theory, with
parametres in ${\frak B}_{\b_{\sig,\a,\zeta+1}}$, satisfied by $\a$, and
saying all the properties of $\a$ we have used so far in this
proof. We can then find $\a'\in
(\b_{\sig,\a,\zeta},\b_{\sig,\a,\zeta+1})$ such that the terms
$\tau_\a$ and $\tau_{\a'}$ coincide, and $m(\a)=m(\a')$,
$n(\a)=n(\a')$, $\langle\xi^{\a'}_1,\dots,\xi^{\a'}_{m(\a')}\rangle =
\langle\xi^\a_1,\dots,\xi^\a_{m(\a)}\rangle$. For every $\xi\leq\sig$,
$h$ maps $M_{\b_{\sig,\a,\xi}}$ into $EM(\b_{\sig,\a,\xi},\Phi')$,
because $\b_{\sig,\a,\xi}\in E^*$.

Now compare the types of $h(a_\a)$ and $h(a_{\a'})$ on
$h(N_{\a'})\subset EM(\a',\Phi)$.

The first one does not $\mu$-split
by monotonicity and the choice of $\zeta$, whereas the second one
$\mu$-splits by the construction, as $p$ $\mu$-splits over
$M_\zeta$. This contradicts the fact that the two types are the same
by the way $\a'$ was chosen.
\qed$_{\ref{splitting}}$

\section{Building the right kind of limits}

We build here from the bottom up the right kind of limit, in
order to approach the construction of models with strong saturation.

\begin{hypothesis}
\end{hypothesis}
\begin{description}
\item[(a)] $LS(\GK )\leq \mu$,
\item[(b)] On $\mu$ we have the consequences of \S 1 and \S 2, namely
density of amalgamation (\ref{density}) and non $\mu$-splitting
(\ref{splitting}).
\item[(c)] Categoricity in $\l > \mu$, $\l \geq \beth_{(2^{LS(\GK
)})^+}$, or at
least some consequences of this.
\end{description}

\subsection{Good extensions. Towers for Limits.}

\begin{definition}
For $\a < \mu^+$, let
\begin{description}
\item[(a)] 
$\GK _{\mu,\a} = \left \{ (\vec{M},\vec{a}) \left |
\begin{tabular}{l}
$\vec{M} = \langle M_i |
i<\a \rangle$ is $\subm$-increasing (not necessarily\\
continuous), $\vec{a} = \langle a_i | i+1<\a \rangle$, $a_i \in
M_{i+1}\setminus M_i$,\\
$M_i\in \GK_\mu$
\end{tabular}\right .\right \}$,


\item[(b)] $\GK ^{am}_{\mu,\a} = \{ (\vec{M},\vec{a}) \in \GK
_{\mu,\a}|$ each $M_i$
is an amalgamation base$\}$,

\item[(c)] $\GK ^{\theta}_{\mu,\a} = \{ (\vec{M},\vec{a}) \in \GK
_{\mu,\a}|$ each $M_i$
is a $(\mu,\theta)$-limit$\}$,

\item[(d)] $$\GK ^*_{\mu,\a} = \bigcup _{\theta \in \mu^+ \cap {\rm Reg}}
\GK ^{\theta}_{\mu,\a},$$
where {\rm Reg}  denotes the class of regular cardinals,
\item[(e)]
${}^+\GK ^{\theta}_{\mu,\a} = \left \{ (\vec{M},\vec{a},\vec{N}) \left |
\begin{tabular}{l}
$(\vec{M},\vec{a}) \in \GK ^{\theta}_{\mu,\a}$, $\vec{N}
= \langle N_i | i+1<\a \rangle$, $N_i \subm M_i$,\\
$N_i$ an amalgamation base in $\GK _\mu$,\\
$M_i$ universal over $N_i$,\\
$\tp (a_i,M_i,M_{i+1})$ does not $\mu$-split over $N_i$
\end{tabular}\right .\right \}$,


\item[(f)] $${}^+\GK ^*_{\mu,\a} = \bigcup _{\theta \in \mu^+ \cap {\rm Reg}}
{}^+\GK ^{\theta}_{\mu,\a}.$$
\end{description}
\end{definition}

\begin{remark}
\end{remark}
\begin{description}
\item[(1)] It is worth noting that, unlike what was done in other
treatments of the subject (see, for example, \cite{Sh:87a},
\cite{Sh:87b}, \cite{Sh:88}, \cite{MaSh:285}, \cite{KlSh:362},
\cite{Sh:394}, \cite{Sh:472} and \cite{Sh:576}), here from now on we
mainly deal with {\bf towers} of models. Objects akin to the towers
defined here were also used in
\cite[\S8-\S10]{Sh:576} in a different context: there full
amalgamation is obtained, but for very few cardinals (only 3 of
them!)... here, we only have amalgamation for {\bf dense} families of
models, but for many more cardinals. We aim at obtaining in subsequent
papers a full description of the categoricity spectrum; in that
respect, amalgamation is a central feature. On the other hand, in
\cite[\S8-\S10]{Sh:576}, the construction is used in order to get the
{\bf non-forking amalgamation}, which is far down the road yet in our
situation.
\item[(2)] What is the point of the definition of ${}^+\GK
^*_{\mu,\a}$? The idea is that we intend to have a parallel to
`the stationarisation
of $\tp (a_i,M_i,M_{i+1}) \in {\cal S}(M_i')$, whenever $M_i \subm
M_i' \in \GK ^{am}_\mu$'. We now turn to defining three orders on the
previously defined classes of towers of models. With these orderings
we intend to capture strong enough notions of limit.
\item[(3)] {\bf Continuity} is not demanded in the definitions
above. One of the major aims is to show that the continuous
towers are dense.
\end{description}

\bigskip

\begin{definition} For $\ell = 1,2$,
\begin{description}
\item[1)] for $(\vec{M}^\ell,\vec{a}^\ell)\in \GK_{\mu,\a}$, let
$$(\vec{M}^1,\vec{a}^1)\suba (\vec{M}^2,\vec{a}^2)$$
mean $\vec{a}^1 = \vec{a}^2$, and for all $i<\a$, $M^1_i\subm M^2_i$,
\item[2)] for $(\vec{M}^\ell,\vec{a}^\ell)\in \GK^{am}_{\mu,\a}$, let
$$(\vec{M}^1,\vec{a}^1)\subb (\vec{M}^2,\vec{a}^2)$$ mean $\vec{a}^1 =
\vec{a}^2$, and for all $i<\a$, $M^1_i=M^2_i$ \underline{or}
$M^1_i\subm M^2_i$, and moreover $M^2_i$ is universal over $M^1_i$,
\item[3)] for $(\vec{M}^\ell,\vec{a}^\ell,\vec{N}^\ell)\in
{}^+\GK^*_{\mu,\a}$, let
$$(\vec{M}^1,\vec{a}^1,\vec{N}^1)\subc
(\vec{M}^2,\vec{a}^2,\vec{N}^2)$$ mean $\vec{a}^1 = \vec{a}^2$,
$\vec{N}^1 = \vec{N}^2$ and for all $i<\a$, $M^1_i=M^2_i$
\underline{or} $M^1_i\subm M^2_i$,
$M^2_i$ is universal over $M^1_i$ (in $\GK _\mu$) and $\tp
(a^1_i,M^2_i,M^2_{i+1})$ does not $\mu$-split over $N^1_i$,
\item[(4)] In all these cases, we
say `strictly' and write `$<^x_{\mu,\a}$, for $x = a$, $b$ or $c$'
if $\bigwedge\limits _iM^1_i\not= M^2_i$).
\end{description}
\end{definition}

\bigskip

\noindent
We have the following facts.

\begin{fact}\label{inclusion}
$\GK _{\mu,\a} \supset \GK ^{am}_{\mu,\a} \supset \GK ^*_{\mu,\a}$.
\end{fact}

\Proof The second inclusion is due to Fact \ref{amalg3}.
\qed$_{\ref{inclusion}}$

\bigskip

\begin{fact}\label{basic}
\begin{description}
\item[1)] $\suba$ is a partial order,
\item[2)] $\subb$ is a partial order,
\item[3)] $\subb \subset \suba$,
\item[4)] If $\langle (\vec{M}^\zeta,\vec{a}^\zeta) | \zeta <\del
\rangle$ is a $\subb$-increasing sequence of members of $\GK
^{am}_{\mu,\a}$, $\del$ is a limit ordinal $<\mu ^+$, and
$(\vec{M},\vec{a}) = \big( \langle
\bigcup_{\zeta<\del}M^\zeta_i|i<\a \rangle,\vec{a}^\zeta \big)$,
$\then$
\begin{enumerate}
\item[(a)] $(\vec{M},\vec{a})\in \GK^{am}_{\mu,\a}$,
\item[(b)] $(\vec{M},\vec{a})$ is the least upper bound of $\langle
(\vec{M}^\zeta,\vec{a}^\zeta) | \zeta <\del \rangle$ (both in
$\suba$ and $\subb$).
\end{enumerate}
\end{description}
\end{fact}
\qed$_{\ref{basic}}$

\begin{remark}
\end{remark}
Part (4) of \ref{basic} explains why we need to have $\subb$ in
addition to $\suba$.

\bigskip
\noindent
This lists the main basic properties of $\suba$, $\subb$ and
$\subc$. It is worth mentioning here that \ref{basic}~(4) has several
uses in what will come next.

\bigskip
\begin{fact}\label{nonempty}
$\GK ^*_{\mu,\a}$, ${}^+\GK ^*_{\mu,\a}$ are both non-empty.
\end{fact}

\Proof We construct the sequences `from the bottom up'. Choose (by
induction on $i<\a$) $M_i\in \GK _\mu$,
$\subm$-increasing continuous, such that
\begin{description}
\item $M_0$ is $(\mu,\omega)$-limit,
\item $M_{i+1}$ is $(\mu,\omega)$-limit and universal over $M_i$,
\item for $i$ limit, $M_i$ is chosen by continuity.
\end{description}
Choose $a_i\in M_{i+1}\setminus M_i$, and choose $N_i$ by using
\ref{universal} and \ref{splitting}. It is easy to see that the resulting sequence of
`double towers' $\langle (M_i,a_i,N_i) | i<\a \rangle$ belongs to
${}^+\GK ^*_{\mu,\a}$, and the corresponding $\langle (M_i,a_i) | i<\a
\rangle$ to $\GK^*_{\mu,\a}$.
\qed$_{\ref{nonempty}}$

\bigskip
We now get a weak form of disjoint amalgamation.

\begin{theorem}\label{open}
If $M_0$ is $(\mu,\theta)$-limit, $M_0 \subm M_\ell$, $M_\ell \in \GK
_\mu$, for $\ell =1, 2$ and $b\in M_1$ \then we can find $M_3$, with
$M_1\subm M_3 \in \GK _\mu$ and a $\subm$-embedding $h$ of $M_2$ into
$M_3$ such that $b\notin h``(M_2)$.
\end{theorem}

\Proof Suppose not. Then fix $M_0$, $M_1$, $M_2$ as in the statement, and
for $i<\mu^+$, find $N_i \in \GK_\mu$, $\subm$-increasing continuous,
and additionally, also find $N^0_i$, $N^1_i$, $N^2_i$ whenever $\cf
(i) = \theta$, such that every $N_i$ is an amalgamation base,
$N_{i+1}$ is universal over $N_i$, and
$$\cf(i) = \theta \implies N_i = N^0_i \subm N^\ell _i
\subm N_{i+1}, \ell = 1, 2,$$
and $(N^0_i,N^1_i,N^2_i,b_i)\approx
(M_0,M_1,M_2,b)$, for some $b_i\in N^1_i$.


Without loss of generality, $$N := \bigcup_{i<\mu^+}N_i\subm
EM(\mu^+,\Phi).$$ Let $E\subset \mu^+$ be a club thin enough so that,
in particular,
$$\del\in E \implies N\cap EM(\del,\Phi) = N_\del.$$
Let also $b_i = \tau_i(\a_{i,0},\dots ,\a_{i,n_i -1})$, with $\a_{i,m_i
-1}<i\leq \a_{i,m_i}$, and $$\bigwedge_{\ell <n_i -1}\a_{i,\ell} <
\a_{i,\ell +1}.$$ Now choose $\del_0\in E$, with $\cf
(\del_0)=\theta$, $\del_0<\del_1\in E$. Let $h$ be the
$\subm$-mapping, with $\Dom h = EM(\del_1,\Phi)$, induced by

\medskip
$\quad \quad \quad \quad \quad \quad j \mapsto \left \{
\begin{tabular}{ll}
$j$&\ife $j<\del_0$\\
$\del_1 + j$&\ife $\del_0 \leq j < \del_1$.

\end{tabular}\right .$

\medskip
\noindent
On $(N^0_{\del_0},N^1_{\del_0},N^2_{\del_0},b_{\del_0})$, we get
precisely the required embedding, and this contradicts the
assumption of its non-existence.
\qed$_{\ref{open}}$

\bigskip

\begin{fact}\label{existgood}{\bf (Existence of good extensions)}
\begin{description}
\item[1)] If $(\vec{M},\vec{a})\in \GK ^*_{\mu,\a}$, and $\theta \in
\mu^+\cap {\rm Reg}$, \then there is
$(\vec{M'},\vec{a'})$ with $(\vec{M},\vec{a})<^b_{\mu,\a}
(\vec{M'},\vec{a'})\in \GK^\theta_{\mu,\a} (\subset \GK ^*_{\mu,\a})$,
where $(\vec{M},\vec{a})<^b_{\mu,\a} (\vec{M'},\vec{a'})$ means
$(\vec{M},\vec{a})\leq^b_{\mu,\a}(\vec{M'},\vec{a'})$ and
$\bigwedge\limits_{\b<\a}[M_\b\not= M_\b']$.
\item[2)] Similarly for ${}^+\GK ^*_{\mu,\a}$, ${}^+\GK
^\theta_{\mu,\a}$, and $\subc$.
\end{description}
\end{fact}

\Proof
\begin{description}
\item 1) Start by observing that given any $M\in \GK^{am}_\mu$, there
is $M'\in \GK_\mu$ universal over $M$ which is actually a
$(\mu,\th)$-limit over $M$: just apply $\th$ many times
\ref{universal} (Existence of Universal Extensions). We still need to
ensure that we get the `weak disjoint amalgamation property', namely
$a_i\notin M_i$'. Theorem~\ref{open} exactly provides this.
\item 2) Like 1), together with the existence of stationarisation of
types and the locality of non-$\mu$-splitting (\ref{splitting}).
\end{description}
\qed$_{\ref{existgood}}$

\bigskip
We now get even more about the least upper bounds for the order
$\subc$.
\bigskip

\begin{fact}\label{basicc}
\begin{description}
\item[1)] $\subc$ is a partial order,
\item[2)] If $\langle (\vec{M}^\zeta,\vec{a}^\zeta,\vec{N}^\zeta) |
\zeta < \del\rangle$ is a $\subc$-increasing sequence of members of
${}^+\GK ^*_{\mu,\a}$, $\del$ is a limit $<\mu^+$, and
$(\vec{M},\vec{a})$ is as in \ref{basic}~(4), \then

\begin{enumerate}
\item[(a)] $(\vec{M},\vec{a})\in \GK^*_{\mu,\a}$,
$(\vec{M},\vec{a},\vec{N})\in {}^+\GK^*_{\mu,\a}$.
\item[(b)] $(\vec{M},\vec{a})$ is the l.u.b. of $\langle
(\vec{M}^\zeta,\vec{a}^\zeta) | \zeta <\del \rangle$ (both in $\suba$
and $\subb$), and
\item[(c)] $(\vec{M},\vec{a},\vec{N})$ is also a $\subc$-l.u.b. of $\langle
(\vec{M}^\zeta,\vec{a}^\zeta,\vec{N}^\zeta) | \zeta <\del \rangle$,
where $\vec{N} = \vec{N}^\zeta$, for any $\z$ (remember they are all equal).
\end{enumerate}
\end{description}
\end{fact}

\Proof
\begin{description}
\item[1)] Trivial,
\item[2)] If the
conclusion were not to hold, then we would fall into
`possibility (a)' of the proof of \ref{splitting}, namely: if $\langle
M_i | i\leq\sig\rangle$ is $\subm$-increasing and continuous, and for
all $i\leq\sig$, $M_i\in\GK^{am}_\mu$, $M_{i+1}$ is universal over
$M_i$, $p\in {\cal S}(M_\del)$, and $p\rest M_i$ does not $\mu$-split
over $M_0$. But then, using \ref{splitting}, we have that $p$ does
not $\mu$-split over $M_0$.

\end{description}
\qed$_{\ref{basicc}}$

\bigskip

\begin{definition}
\begin{description}
\item[1)] $(\vec{M},\vec{a})\in \GK^*_{\mu,\a}$ is {\bf reduced} if
$$(\vec{M},\vec{a})\subb (\vec{M'},\vec{a'})\implies
\bigwedge_{i<\a} \Big[ M'_i \cap \bigcup_{j<\a}M_j = M_i\Big] .$$
\item[2)] $(\vec{M},\vec{a},\vec{N})\in {}^+\GK^*_{\mu,\a}$ is {\bf
reduced} if $$(\vec{M},\vec{a},\vec{N})\subc (\vec{M'},\vec{a},\vec{N})
\implies \bigwedge_{i<\a} \Big[ M'_i \cap \bigcup_{j<\a}M_j = M_i\Big]
.$$
\end{description}
\end{definition}

\begin{remark}
\end{remark}
Equivalently, when defining $(\vec{M},\vec{a})$ is reduced, we could
have used $\suba$ instead of $\subb$: just notice that $x\suba x'\subb
x'' \implies x\subb x''$ and for all $x$ in the appropriate class of
towers there exists $y$ such that $x\subb y$.

\bigskip

\begin{fact}\label{densredtow}{\bf (Density of Reduced Towers)}
\begin{description}
\item[1)] For every $\theta \in \mu^+\cap {\rm Reg}$, for every
$(\vec{M},\vec{a})\in \GK^\theta_{\mu,\a}$, there is a reduced tower
$(\vec{M'},\vec{a'})\in \GK^\theta_{\mu,\a}$ such that
$(\vec{M},\vec{a})\subm (\vec{M'},\vec{a'})$.
\item[2)] Similarly for ${}^+\GK^*_{\mu,\a}$, ${}^+\GK^\theta_{\mu,\a}$.
\end{description}
\end{fact}

\Proof \begin{description}
\item[(1)] Let $(\vec{M},\vec{a})\in\GK^\th_{\mu,\a}$. If the
conclusion fails, then we can find
$(\vec{M}^i,\vec{a})\in\GK^\th_{\mu,\a}$, $\subb$-increasing
continuous for $i<\mu^+$, such that $(\vec{M}^{i+1},\vec{a})$
witnesses that $(\vec{M}^i,\vec{a})$ is not reduced. Now the set
$$E= \Big\{ \del<\mu^+ \Big| i<\a \implies \Big(
\bigcup_{\z<\mu^+}M^\z_i\Big) \cap \Big(
\bigcup_{j<\a}M^\del_j\Big) = M^\del_i \Big\}$$ is a club of
$\mu^+$. For $a\in \bigcup\limits_{\stackrel{\z<\mu^+}{i<\a}}M^\z_i$, let
\begin{description}
\item $i(a)=\min \{ i<\a |a\in
\bigcup\limits_{\stackrel{\z<\mu^+}{j<i}}M^\z_j\}$, 
\item $\z(a) = \min \{ \z|a\in \bigcup\limits_{\xi<\z}M^\xi_{i(a)}\}$.
\end{description}
So, $E' = \{\del |a\in M^\del_i$, $i<\a \implies \z(a)<\del\}$ is a
club and $E'\subset E$.  Choose $\del^*\in E'$: this violates that the
conclusion fails. [Why? $(\vec{M}^{\del^*},\vec{a})$ is reduced. To see this,
just let $(\vec{M}^{\del^*},\vec{a})\subb (\vec{M}',\vec{a})$. We only
have to check that for every $i<\a$, $M'_i\cap
\bigcup\limits_{j<\a}M_j^{\del^*}\subset M^{\del^*}_i$. So let $a\in M'_i\cap
\bigcup\limits_{j<\a}M_j^{\del^*}$. Since $\del^*\in E'$,
$\z(\a)<\del^*$, hence $a\in
\bigcup\limits_{\xi<\z(\a)}M_{i(a)}^\xi$. But since $\del^*\in E$, we
have that $\Big( \bigcup\limits_{\z<\mu^+}M^\z_i\Big) \cap \Big(
\bigcup\limits_{j<\a}M^{\del^*}_j\Big) = M^{\del^*}_i$. This implies that $a\in
M^{\del^*}_i$.]

\item[(2)] Clearly similar.
\end{description}
\qed$_{\ref{densredtow}}$

\begin{fact}\label{basicred}
\begin{description}
\item[1)] In \ref{basic}~(4), if $\del = \sup \{\zeta<\del |
(\vec{M}^\zeta,\vec{a}^\zeta)$ is reduced$\}$, \then
$(\vec{M},\vec{a})$ is reduced. (In fact, it is enough to have
$(\vec{M},\vec{a})\in \GK^{am}_{\mu,\a}$, and
$(\vec{M}^\zeta,\vec{a}^\zeta)\in \GK^{am}_{\mu,\a}$ is
$\suba$-increasing.)
\item[2)] In \ref{basicc}~(2), if $\del = \sup \{\zeta<\del |
(\vec{M}^\zeta,\vec{a}^\zeta,\vec{N}^\zeta)$ is reduced$\}$, \then
$(\vec{M},\vec{a},\vec{N}^\zeta)$ is reduced.
\end{description}
\end{fact}

\Proof Clear from the definition of `reduced'. \qed$_{\ref{basicred}}$

\bigskip
\begin{theorem}\label{reducont}
\begin{description}
\item[1)] If $(\vec{M},\vec{a})\in \GK^*_{\mu,\a}$ is reduced, \then
$\vec{M}$ is $\subm$-increasing {\it and\/} continuous.
\item[2)] If $(\vec{M},\vec{a},\vec{N})\in {}^+\GK^*_{\mu,\a}$ is
reduced, \then $\vec{M}$ is $\subm$-increasing {\it and\/} continuous.
\end{description}
\end{theorem}

\Proof We prove by induction on $\del<\a$ limit ordinal that if
$(\vec{M},\vec{a})\in \GK^*_{\mu,\a}$ is reduced, then $M_\del =
\bigcup\limits_{i<\del}M_i$. Assume then failure for $\del$: there
exists some $b\in M_\del\setminus \bigcup\limits_{i<\del}M_i$. We can
find $(\vec{M}^\zeta,\vec{a})$ {\it reduced\/}, for $\zeta\leq \del$,
a $\subb$-strictly increasing continuous chain of towers, such that
$(\vec{M}^0,\vec{a}) = (\vec{M},\vec{a})$ (it exists by
\ref{basic}~(4), and because `reduced towers are dense and closed under
limit', \ref{densredtow}). Now consider the diagonal sequence $\langle
M^\zeta_\zeta |
\zeta\leq\del\rangle$. It is $\subm$-increasing continuous, its
members are in $\GK _\mu$, each $M^\zeta_\zeta$ is
$(\mu,\th_\z)$-limit, for some $\th_\z\in \mu^+\cap {\rm Reg}$, and
$M^{\z+1}_{\z+1}$ is universal over $M^\z_\z$. Also,
$M^\del_\del\subset M^\del_{\del+1}$, and $b\in M^\del_{\del+1}$. So,
by the main result on non-$\mu$-splitting (Theorem~\ref{splitting}),
for some $\xi<\del$, $\tp
\big( b,\bigcup\limits_{\xi<\del}M^\xi_\xi,M^\del_\del\big)$ does not
$\mu$-split over
$M^\xi_\xi$. Let now $\frak B$ be $(\mu,\mu^+)$-limit,
$M^\del_\del\subm \frak B$.

\bigskip
We choose by induction on $i\leq \del$ models $N_i$ and functions
$h_i$ such that

\begin{description}
\item[ * ] if $i\leq \xi+1$, \then $N_i = M^\del_i$ and $h_i =$
id$_{N_i}$,
\item[ * ] if $i\in (\xi+1,\del]$, \then $N_i\subm \frak B$ is a
$(\mu,\th_i)$-limit model, and
\begin{enumerate}
\item $\langle N_i|i\leq \del\rangle$ is $\subm$-increasing, continuous,
\item $N_{i+1}$ is universal over $N_i$,
\item $N_{\xi+2} \supset M_\del^\del$,
\item $h_i$ is a $\subm$-embedding of $M_i^\del$ into $N_i$,
\item $h_{i+1}$ maps $M^\del_{i+1}\setminus M^\del_i$ into
$N_{i+1}\setminus N_i$,
\item $\langle h_i|i\leq \del\rangle$ is increasing continuous,
\item $\tp (b,h_i(M_i^\del),{\frak B})$ does not $\mu$-split over $N^*
:= M^\xi_\xi = h_i(M^\xi_\xi)$, for $i\geq \xi$.
\end{enumerate}
\end{description}

\noindent
For $i\leq \xi+1$, this is trivial. For $i\in (\xi+1,\del)$ successor,
by the claim \ref{nonspext} below. For $i\in (\xi+1,\del]$ limit, use
\ref{splitting} for the last clause (\underline{remember}, by the
induction hypothesis, $\langle M_i^\del | i<\del\rangle$ is
continuous, and by definition, $M_i^\del$ is $(\mu,\th)$-limit, hence
an amalgamation base). We also have that $\tp (b,\bigcup\limits
_{i<\del}h_i(M^\del_i),{\frak B})$ does not $\mu$-split over $N^*$. So,
$$h^*=\bigcup\limits_{i<\del}h_i\cup \{(b,b)\}$$
is a `legal' map. For some $N_{\del+1}\subm \frak B$, $(\mu,\w)$-limit
over $N_\del$, we can extend $h^*$ to $h^+\in AUT(N_{\del+1})$. Let
for $i\leq\del$, $M_i^\otimes = (h^+)^{-1}(N_i)$.

\noindent
We then have

\begin{description}
\item[$\bigotimes _b$] $\langle
M_i^\del | i\leq\del\rangle \leq^b_{\mu,\del}\langle M_i^\otimes |
i\leq\del\rangle$
\end{description}

[Why? On the one hand, $i\leq
\xi\implies h^+\supset h_i = id_{M^\del_i} = id_{N_i}$, and thus
$M_i^\otimes = N_i = M^\del_i$. On the other hand, if $i\in
(\xi,\del)$, then $h^+\supset
h^*\supset h_i$ and $h_i$ $\subm$-maps $M^\del_i$ into $N_i$. We thus
have that $M^\del_i\subm M^\otimes_i$. If $i=\del$, then clearly
$M^\del_i\subm M^\otimes_i$.

In the $\GK^*_{\mu,\a}$ case, we still need to show why $i<\del
\implies a_i \notin M^\otimes_i$: if $i\leq \xi$, this is trivial; if
$i>\xi$, as $h_{i+1}$ maps $M^\del_{i+1}\setminus M^\del_i$ into
$N_{i+1}\setminus N_i$, and thus $h_{i+1}(a_i)\notin N_i$, hence
$a_i\notin (h^+)^{-1}(N_i) = M^\otimes_i$.]

\begin{figure}[ht]\hspace*{1.00mm}
\psfig{figure=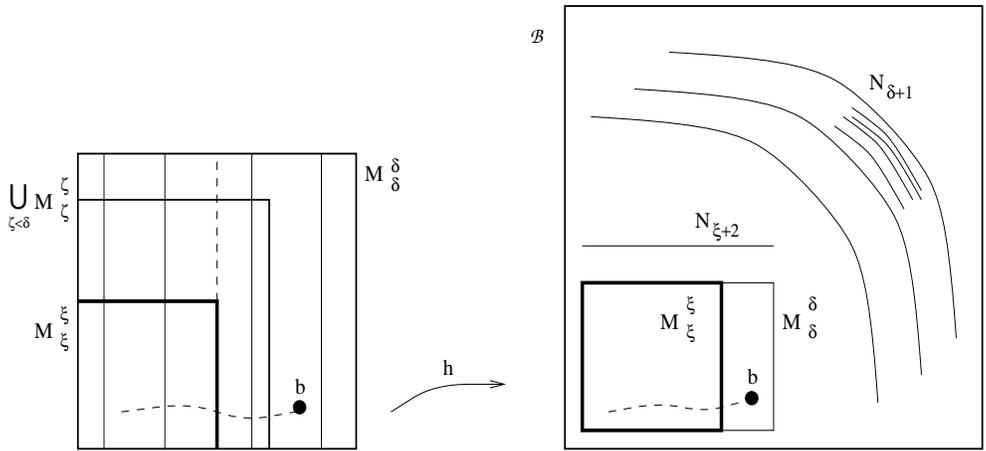,height=6cm,width=13cm}
\caption{ The diagonal and $h$. } 
\end{figure}

\bigskip

So, $\bigotimes _b$ holds, and as earlier we can define $M_i^\otimes$, for
$i\in (\del,\a)$ such that $\langle M_i^\del | i<\a\rangle
\leq^b_{\mu,\del}\langle M_i^\otimes |i<\a\rangle$. But then the place
of $b$ drops: now, $b\in M_\del \subm N_{\xi+1}$, and $h^+(b) = b$, so
$b\in M^\otimes_{\xi+1}$. Contradiction.
This finishes the proof of (1).

\medskip
\noindent
The proof of (2) is similar: We are now in the ${}^+\GK^*_{\mu,\a}$
case, and we need (in addition to what has already proved) to prove
that the non-$\mu$-splitting holds. This is, $\tp
(a_i,M^\otimes_i,M^\otimes_{i+1})$ does not $\mu$-split over
$N^\del_i$. By definition of ${}^+\GK _{\a,\mu}$, we already have that
$\tp (a_i,M^\del_i,M^\del_{i+1})$ does not $\mu$-split over
$N^\del_i$. But $h^+\in AUT(N_{\del+1})$; hence, $\tp
(h^+(a_i),h^+(M^\del_i),h^+(M^\del_{i+1}))$ does not $\mu$-split over
$h^+(N^\del_i)$. But then, pulling back again this type, we have that
$\tp (a_i,M^\otimes_i,M^\otimes_{i+1})$ does not $\mu$-split over
$N^\del_i$.

The proof is finished just like in case (1); we use the fact
\begin{description}
\item[$\bigotimes _c$] $\langle
M_i^\del,N_i^\del | i\leq\del\rangle \leq^c_{\mu,\del}\langle
M_i^\otimes,N_i^\del | i\leq\del\rangle$
\end{description}
the same way as there.
\qed$_{\ref{reducont}}$

\begin{claim}\label{nonspext}
Assume that $M_0$ is $(\mu,\th_0)$-limit over $N_0^*$, $N_0^*\in
\GK^{am}_\mu$, $M_\ell$ is $(\mu,\th_1)$-limit and universal over
$M_0$, for $\ell =1,2$,
$b\in M_1$, $\tp (b,M_0,M_1)$ does not $\mu$-split over $N_0^*$.
\Then, we can find $M_3\in \GK^{am}_\mu$ and $h$ such that $M_1\subm
M_3$, $h$ is a $\subm$-embedding of $M_2$ into $M_3$ over $M_0$, $\tp
(b,h(M_2),M_3)$ does not $\mu$-split over $N_0^*$.
\end{claim}

\Proof We can first find $M_3'$, $h'$ such that $M_2\subm M_3'\in
\GK^{am}_\mu$, $h'$ is a $\subm$-embedding of $M_1$ into $M_3'$ over
$M_0$, and $p'=\tp (h'(b),M_2,M_3')$ does not $\mu$-split over
$N_0^*$. This follows by the definition of type and the existence of
$p'\in {\cal S}(M_2)$, $p'\geq \tp (b,M_0,M_1)$, not
$\mu$-splitting over $N^*_0$: hence, for some $M_3''$, $M_2\subm M_3''\in
\GK^{am}_\mu$, and $b'\in M_3''$ realises $p$. So there are $M_3'$,
$h'$, $M_3'\supset M_3''$ as required: send $b$ to $b'$ via $h'$, and
extend the identity on $M_0$.
\qed$_{\ref{nonspext}}$

\bigskip
\subsection{Toward the uniqueness of limits}

We need a refined concept of type in order to obtain the right kind of
towers later (`full' towers). The following definition specifies this
refinement: in addition to just `describing elements,' like we do when
defining types of various sorts, we look both at the `elements'
themselves and at witnesses of the specific `way they do not
$\mu$-split'. In principle, this provides a tighter description of
the element, since it provides along with it the specific submodel
over which the type does not $\mu$-split.

\bigskip

\begin{definition}\label{st}
For $M$ a $(\mu,\th)$-limit model, let
\begin{description}
\item[1)] $\frak{St}(M) = \{(p,N) | N\subm M$ is $(\mu,\th)$-limit, $M$
universal over $N$, $p\in {\cal S}(M)$ does not $\mu$-split over
$N\}$,
\item[2)] for $(p_\ell,N_\ell)\in \frak{St}(M)$, for $\ell = 1,2$, let
$(p_1,N_1)\approx (p_2,N_2)$ \iff for every $M'$, $M\subm M'\in
\GK^{am}_\mu$, there is $q\in {\cal S}(M')$ extending $p_1$ and $p_2$,
not $\mu$-splitting over $N_1$ or over $N_2$.
\end{description}
\end{definition}

\begin{fact}\label{splitypes}
\begin{description}
\item[1)] $\approx$ is an equivalence relation on $\frak{St}(M)$,
\item[2)] If $M'\in \GK^{am}_\mu$ is universal over $M$ (in
\ref{st}~(2)), the existence of $q$ for this $M'$ suffices,
\item[3)] $\Big| \frak{St}(M)/\approx \Big| \leq \mu$.
\end{description}
\end{fact}

\Proof 
\begin{description}
\item[1)] A diagram chase which we leave to the reader,
\item[2)] By Universality + Preservation by $\subm$-embeddings,
\item[3)] Since by \ref{splitypes}~(2), there is $M'\in \GK^{am}_\mu$
universal over $M$ in which we may check all the instances of
$\approx$-equivalence, we have $\Big|
\frak{St}(M)/\approx \Big| \leq \Big| {\cal S}(M')\Big| \leq \mu$.
\end{description}
\qed$_{\ref{splitypes}}$

\noindent
{\bf Remarks:}
\begin{description}
\item[1)] It is worth noting here that perhaps $\approx$ is the equality. We do
not know yet; but for our purposes, it is OK to use $\approx$.
\item[2)] In the definition of `$(\vec{M},\vec{a},\vec{N})\in
{}^+\GK^*_{\mu,\a}$`, it is just the $\frak{St}(M)$-equivalence class $(\tp
(b_i,M_i,M_{i+1}),N_i)/\approx$ that matters, and not $N_i$ itself.
\end{description}

\bigskip
And now, we can provide a crucial notion for towers (see also \cite{Sh:576}).

\begin{definition}
We say that $(\vec{M},\vec{a},\vec{N})\in {}^+\GK^*_{\mu,\a}$ is {\bf
full} \iff
\begin{description}
\item[(a)] $\mu$ divides $\a$ (if $\mu$ is regular, if it is singular,
$\mu^\w$ divides $\a$).
\item[(b)] if $\b<\a$ and $(p,N^*)\in \frak{St}(M_\b)$, then for some
$i<\b+\mu$, we have that $\Big( \tp
(b_{\b+i},M_{\b+i},M_{\b+i+1}),N_{\b+i}\Big) \approx
(p,N^*)$. (Formally, it is equivalent to the stationarisation of $(p,N^*)$.)
\end{description}
\end{definition}

\bigskip
We are approaching one of our main goals (`uniqueness of limits') with
the following theorems.

\begin{theorem}\label{uniquelimit}
If $(\vec{M},\vec{a},\vec{N})\in {}^+\GK^*_{\mu,\a}$ is full, and $\vec{M}$
is continuous, then $\bigcup\limits_{i<\a}M_i$ is $(\mu,\cf \a)$-limit
over $M_0$.
\end{theorem}

\Proof Let $\langle M_i' | i\leq \a\rangle$ be $\subm$-increasing
continuous, with each model in the tower in $\GK^{am}_\mu$, and such
that each $M_{i+1}'$ is universal over $M_i'$, $M_0'
\stackrel{h_0}{\approx}M_0$. List $|M_i|$ as $\{a_{i,\z} | \z<\mu\}$ and
$|M_i'|$ as $\{a_{i,\z}' | \z<\mu\}$. Let $g:\a\longrightarrow
\bigcup\limits_{i\leq\a}M_i' = M_\a'$ satisfy $g(i)\in M_i'$ and

\bigskip
$\quad \otimes \left [
\begin{tabular}{l}
if $\b<\a$, $b\in M_\b'$, $(p,N)\in \frak{St}(M_\b)$, then\\
\\
$\mu=\otp \left \{ \g \left |
\begin{tabular}{l}
$\b<\g<\b+\mu$\\
$(p,N)\approx \Big( \tp (a_\g,M_\g,M_{\g+1}),N_\g\Big)\restriction M_\b$\\
$g(\g)=b$
\end{tabular}\right .\right \} .$
\end{tabular}\right .$

\bigskip
\noindent
This makes sense: $\| M_\b'\| \leq \mu$ and by \ref{splitypes},
$\Big| \frak{St}(M)/\approx \Big| \leq \mu$.

\begin{figure}[ht]\hspace*{1.00mm}
\psfig{figure=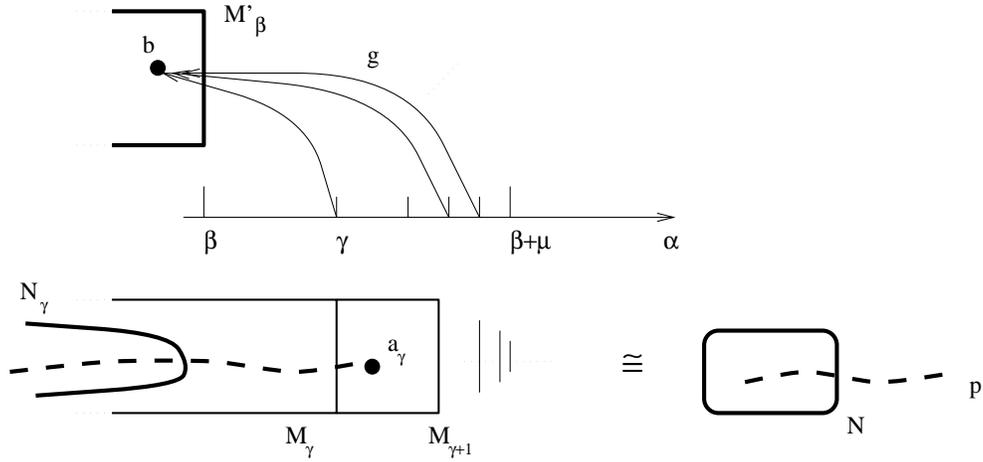,height=6cm,width=13cm}
\caption{ The definition of $g$. } 
\end{figure}

\bigskip
\noindent
We choose by induction on $i<\a$, $h_i$ such that

\begin{description}

\item[(a)] $h_i$ is a $\subm$-embedding of $M_{j^1_i}$ into $M_{j^2_i}'$
\item[(b)] $\langle h_i\rangle _{i<\a}$ is increasing continuous,
\item[(c)] $j^1_i$, $j^2_i$ are increasing continuous,
\item[(d)] $j^1_{i+1}<j^1_i+\mu$, $j^2_{i+1}<j^2_i+\mu$,
\item[(e)] $g(i)\in \Rang (h_{i+1})$.
\end{description}

For $i=0$, this is trivial. for $i$ limit, just take unions. For
\underline{$i+1$}, without loss of generality $g(i)\in M_i'\subm
M_{j^2_i}'$, and we know that $M_i$ is $(\mu,\th_i)$-limit (say
$\langle M_{i,\vare} | \vare\leq \th_i\rangle$ witnesses this). So,
for some $\vare_i$, $\tp (g(i),h_i(M_i),M_{j^2_i}')$ does not
$\mu$-split over $h_i(M_{i,\vare_i})$. So there is $p_i\in {\cal
S}(M_{j^1_i})$ such that $h_i(p_i) = \tp (g(i),h_i(M_i),M_{j^2_i}')$,
and thus for some $\xi\in (0,\mu)$, we have
$$(p_i,M_{i,\vare_i})\approx \Big( \tp
(a_{j^1_i+\xi},M_{j^1_i+\xi},M_{j^1_i+\xi+1}),N_{j^1_i+\xi}\Big ).$$
Let $j^1_{i+1} = j^1_i + \xi +1$, $j^2_{i+1} = j^2_i +1$; there is
$h_{i+1}$ as required.

\bigskip
So, letting $h = \bigcup\limits_{i<\a}h_i$, we have a witness for the
original requirement.

\bigskip
\noindent
For $\mu$ regular, $j^2_i=i\leq j^1_i\in[i,i+\mu)$, for $\mu$
singular, $j^2_i=i$, $j^1_i\in[i,\mu i)$. So if $\mu^\w$ divides $\a$
we catch our tail ($j^2_\del = j^1_\del = \del$), and get the
isomorphism we are looking for ($h_\del:M_\del\longrightarrow M_\del'$
is onto by the bookkeeping $g$).
\qed$_{\ref{uniquelimit}}$

\subsection{Limits via sequences of different lengths}

So far, `the limit' has been proven to be unique, when the sequences
converging to it in the various orderings defined are of same
length. We are striving for more: we want to prove that even if we
approach a model {\bf via sequences of certain different lengths}, the
limit may be proven to be `unique' in a robust enough sense, by using
a rectangle of models which will be $(\mu,\th_\ell)$-limits over
$M_{0,0}$, for $\ell=1,2$, by the two sides.

\begin{definition} For $u$ an interval, and $\gu$ a set of intervals, let
\begin{description}
\item ${}^+\GK ^*_{\mu,u} = \left \{ (\vec{M},\vec{a},\vec{N}) \left |
\begin{tabular}{l}
$\vec{M} = \langle M_i | i\in u \rangle$ is $\subm$-increasing (not\\
necessarily continuous), $\vec{a} = \langle a_i | i\in u \rangle$,\\
$a_i \in M_{i+1}\setminus M_i$, $\vec{N} = \langle N_i | i\in u
\rangle$, $N_i \subm M_i$,\\
$N_i$ an amalgamation base in $\GK _\mu$,\\
$M_i$ universal over $N_i$,\\
$\tp(a_i,M_i,M_{i+1})$ does not $\mu$-split over $N_i$
\end{tabular}\right .\right \}$,

\item ${}^+\GK ^*_{\mu,\gu} = \bigcup\limits_{u\subset \gu}{}^+\GK ^*_{\mu,u}$.
\end{description}
\end{definition}

\noindent
The right way to think about these classes is by immediate analogy to
the original ${}^+\GK ^*_{\mu,\a}$ classes. As there, it is natural to
expect to have a
$\subc$ relation.

\begin{definition}
For $(\vec{M}^\ell,\vec{a}^\ell,\vec{N}^\ell)\in {}^+\GK
^*_{\mu, \gu_\ell}$, $\ell = 1,2$, let
$$(\vec{M}^1,\vec{a}^1,\vec{N}^1)\leq^c(\vec{M}^2,\vec{a}^2,\vec{N}^2)$$
mean that $(\vec{M}^1,\vec{a}^1,\vec{N}^1)\subc (\vec{M}^2\rest
\gu_1,\vec{a}^2\rest \gu_1,\vec{N}^2\rest \gu_1)$.
\end{definition}

Now, as before for the definitions of towers, we have the following
basic facts about the new `scattered' towers.

\begin{fact}\label{fullbasic}
\begin{description}
\item[1)] Let $\gu_1\subset \gu_2\subset ORD$. Then,
\begin{enumerate}
\item[(a)] If $(\vec{M},\vec{a},\vec{N})\in {}^+\GK ^*_{\mu,\gu_2}$,
then $(\vec{M},\vec{a},\vec{N})\rest \gu_1 \in {}^+\GK
^*_{\mu,\gu_1}$,
\item[(b)] If $(\vec{M},\vec{a},\vec{N})\in {}^+\GK ^*_{\mu,\gu_1}$,
then there is $(\vec{M}',\vec{a}',\vec{N}')\in {}^+\GK ^*_{\mu,\gu_2}$
such that $(\vec{M}',\vec{a}',\vec{N}')\rest \gu_1\in
{}^+\GK ^*_{\mu,\gu_2}$.
\end{enumerate}
\item[2)] Let $\langle \gu_\vare | \vare <\vare (*)\rangle$ be an
increasing sequence of sets of ordinals such that $|\gu_\vare|\leq \mu$
then the parallel of \ref{basicc} for limits.
\item[3)] If $\gu_1\subset \gu_2$,
$(\vec{M}^\ell,\vec{a}^\ell,\vec{N}^\ell)\in {}^+\GK
^*_{\mu,\gu_\ell}$, and $(\vec{M}^2,\vec{a}^2,\vec{N}^2)\rest \gu_1 \leq ^c
(\vec{M}^1,\vec{a}^1,\vec{N}^1)$, \then we can find
$(\vec{M}^3,\vec{a}^3,\vec{N}^3)\in {}^+\GK ^*_{\mu,\gu_1}$ such that
$(\vec{M}^1,\vec{a}^1,\vec{N}^1)\leq ^c
(\vec{M}^3,\vec{a}^3,\vec{N}^3)$
and $(\vec{M}^2,\vec{a}^2,\vec{N}^2)\leq ^c (\vec{M}^3,\vec{a}^3,\vec{N}^3)$.
\item[4)] If $\gu\subset \gu_1\cap \gu_2$,
$(\vec{M}^\ell,\vec{a}^\ell,\vec{N}^\ell)\in {}^+\GK
^*_{\mu,\gu_\ell}$, for $\ell = 1,2$, $(\vec{M}^2,\vec{a}^2,\vec{N}^2)\rest
\gu_1 \leq ^c(\vec{M}^1,\vec{a}^1,\vec{N}^1)$, then
$(\vec{M}^2,\vec{a}^2,\vec{N}^2)\rest \gu \leq
^c(\vec{M}^1,\vec{a}^1,\vec{N}^1)\rest \gu$.
\end{description}
\end{fact}

\noindent
We also have that
\begin{fact}
In \ref{fullbasic}, if each one of the towers is reduced, then so are
the limits.
\end{fact}

\bigskip

\begin{construction}\label{construction}
\end{construction}
Fix $\z$ and let $\a<\mu^+$ be divisible by $\mu^\w$ (ordinal
exponentiation). Let now ${\frak U}_\vare = \bigcup \big\{
[\b\mu\z,\b\mu\z +\mu\vare) | \b<\a\}$, for each $\vare\leq \z$. Now
we can define the class ${}^+\GK^*_{\mu,\gu_\vare}$ just like
${}^+\GK^*_{\mu,\a}$, but now using $\gu_\vare$ as our set of indices,
instead of $\a$. The point is that now we want to play with changing
$\gu_\vare$ in various ways.

\bigskip
\noindent
We choose $(\vec{M}^\vare,\vec{a}^\vare,\vec{N}^\vare)\in {}^+\GK
^*_{\mu,\gu_\vare}$, $<^c$-increasing (naturally), each one of these towers
reduced. In successor stages, in the new intervals which have length
$\mu$, we put representatives to all types.

\bigskip
In the end, we get $\vec{M} = \langle M_i | i<
(\b\mu\z)\a\rangle$. So,

\begin{description}
\item[$\otimes _1$] each $M_i$ is $(\mu, \cf \z)$-limit (by
$<^c$),
\item[$\otimes _2$] every $(p,N)\in \frak{St}(M_i)/\approx$, up to
equivalence (by~\ref{splitting} --- we dealt with it so it appears in
$j\in [i,i+\b\mu\z)$),
\item[$\otimes _3$] if there exists $\vare$ such that $i\in \gu_\vare$
and $i=\sup (\gu_\vare \cap i)$ then $M_i = \bigcup\limits_{j<i}M_j$.
\end{description}

So, reformulating `full $\implies$ limit,' we get a similar claim for
$\gu$ instead of $\mu$ or $\mu ^\omega$.

If $\a=\a^*+1$ is a large enough ordinal, then we can find $\gv$
such that $\vec{M}\rest \gv$ is full, $\a^*=\sup \gv\cap
\a^*$ and $\a^*\in \gv$. So, by $\otimes _3$, we have

\begin{description}
\item[$\bigotimes$] $M_{\a^*} = \bigcup\limits_{\b\in\gv}M_\b$ is
$(\mu,\cf \a^*)$-limit over $M_0$.
\end{description}

\noindent
But $M_\a$ is $(\mu,\cf \z)$-limit. We can arrange $\cf \a^*$, $\cf
\z$ to be any regular $<\mu^+$.

\begin{conclusion}
\end{conclusion}
If $M_\ell$ is $(\mu,\th_\ell)$-limit ($\ell = 1,2$), then $M_1
\approx M_2$.

\bigskip
\noindent
So, to speak about {\bf `the $\mu$-limit model'} now makes sense.

\bigskip
\noindent
{\bf Comments:} A nicer construction may be obtained if we set

\bigskip
$\gu_\vare = \left \{ \mu^\w\z + j \left |
\begin{tabular}{ll}
$j<\mu^\w\z$,&$j=0 \mod 3$\\
	     &or $j=1 \mod 3$\\
	     &or $j<\mu^\w\vare$
\end{tabular}\right .\right \}$.

\bigskip
This way, the first set $\gu_0$ contains at least all the
ordinals which are $0$ or $1\mod 3$, hence there is no problem with
the limit.

\begin{theorem}\label{satur}
If $M_i$ is $(\mu,\th_\ell)$-limit over $M$, for $\ell =1,2$, then
$$M_1 \approx _M M_2.$$
\end{theorem}

\Proof The same proof as for Theorem \ref{uniquelimit} works, although
naturally it has to be adapted to our `scattered tower'
situation. Without loss of generality, both $0$ and $\a^*\in \gu$. We
define $(\vec{M}^\z,\vec{a}^\z,\vec{N}^\z)\in {}^+\GK
^*_{\mu,\gu_\z}$, $M^0_0= M$. So, we have that $M^\z_{\a^*}$ is a
$(\mu,\th_{{\rm cf}\a^*})$-limit over $M$, and also $(\mu,\cf
\z)$-limit over $M$.
\qed$_{\ref{satur}}$

\bigskip
With this, we can by now conclude that {\bf saturated models exist} in
a strong enough sense. We may take as our definition of a `saturated
model over $M$' in a cardinal $\mu$ the (by now unique because of
Theorem~\ref{satur}) $(\mu,\th)$-limit over $M$, for an arbitrary
$\th$.

\bigskip
This paves the way towards a full study of the categoricity spectrum
for abstract elementary classes without maximal elements. We intend to
continue developing this theory in that direction, by studying the
type theory for our context, non forking amalgamation, and the true
role of saturation. But this will be the material of forthcoming papers.

\bibliographystyle{literal-plain}

\bibliography{lista,listb,listf,listx}
\end{document}